# *On sandwich theorems for univalent meromorphic functions involving integral operator*

Waggas Galib Atshan and Khudair O. Hussain

Department of mathematics

College of Computer Science and Information Technology

University of Al-Qadisiya, Diwaniya, Iraq

waggashnd@gmail.com , khudair_o.hussain@yahoo.com

**Abstract**

*The main purpose of this paper is to derive some subordination and superordination results involving certain of integral operator for meromorphic univalent functions in the punctured open unit disk. Several sandwich-type results are also obtained.*

**Keywords:** meromorphic functions, subordination, superordination, sandwich Theory, integral operator.

**AMS Subject Classifications: 30C45**

## Introduction

Let $\Sigma$ denote the class of functions of the form:

$$f(z) = z^{-1} + \sum_{k=1}^{\infty} a_k z^k, \qquad (1.1)$$

which are analytic in the punctured open unit disk $U^* = \{z: z \in \mathbb{C},\ 0 < |z| < 1\}$.

Let $H$ be the linear space of all analytic functions in $U$. For a positive integer number $n$ and $a \in \mathbb{C}$, we let

$$H[a, n] = \{f \in H: f(z) = a + a_n z^n + a_{n+1} z^{n+1} + \cdots\}.$$

For two functions $f$ and $g$ analytic in $U$, we say that the function $g$ is subordinate to $f$ in $U$ and write $g(z) \prec f(z)$, if there exists a Schwarz function $\omega$, which is analytic in $U$ with $\omega(0) = 0$ and $|\omega(z)| < 1$ $(z \in U)$, such that $g(z) = f(\omega(z))$, $(z \in U)$. Indeed, it is known that

$$g(z) \prec f(z) \Rightarrow g(0) = f(0) \text{ and } g(U) \subset f(U).$$

Furthermore, if the function $f$ is univalent in $U$, we have the following equivalence relationship (cf., e.g. [3], [7] and [8])

$$g(z) \prec f(z) \Leftrightarrow g(0) = f(0) \text{ and } g(U) \subset f(U), \qquad (z \in U).$$

Supposing that $p$ and $h$ are two analytic functions in $U$, let $\emptyset(r, s, t; z) = \mathbb{C}^3 \times U \to \mathbb{C}$,

if $p$ and $\emptyset(p(z), zp'(z), z^2 p''(z); z)$ are univalent functions in $U$ and $p$ satisfies the second order superordination

$$h(z) \prec \emptyset(p(z), zp'(z), z^2 p''(z); z), \qquad (z \in U). \qquad (1.2)$$

Then $p$ is called a solution of the differential superordination (1.2). (if $f$ subordinate to $g$, then $g$ is superordinate to $f$).

An analytic function $q$ is called a subordinate of the differential superordination if $q \prec p$ for all $p$ satisfying (1.2). A univalent subordinate $\tilde{q}$ that satisfies $q \prec \tilde{q}$ for all subordinates $q$ of (1.2) is said to be the best subordinate. Recently Miller and Mocnuu [8]obtained sufficient conditions on the functions $h, p$ and $\emptyset$ for which the following implication holds

$$h(z) \prec \emptyset(p(z), zp'(z), z^2 p''(z); z) \Longrightarrow q(z) \prec p(z) \qquad (z \in U) \qquad (1.3)$$

Using the results, Bulboaca[2]considered certain classes of first order differential superordinations as well as superordination preserving integral operator [3].Ali et al.[1], have used the results of Bulboaca [2]to obtain sufficient conditions for normalized analytic functions to satisfy

$$q_1(z) \prec \frac{zf'(z)}{f(z)} \prec q_2(z),$$

where $q_1$ and $q_2$ are given univalent functions in $U$ with $q_1(0) = q_2(0) = 1$. Also, Tuneski [11] obtained a sufficient condition for starlikeness of $f$ in terms of the quantity

$$\frac{f''(z)f(z)}{(f'(z))^2}.$$

Recently, Shanmugam *et al.* [9,10] and Goyal *et al.*[4] also obtained sandwich results for certain classes of analytic functions.

A.Y. Lashin [5] introduced and investigated the integral operator

$$\mathcal{P}_\beta^\alpha : \Sigma \longrightarrow \Sigma$$

which is defined as follows:

$$\mathcal{P}_\beta^\alpha f(z) = \frac{\beta^\alpha}{\Gamma(\alpha)} \frac{1}{z^{\beta+1}} \int_0^z t^\beta \left(\log \frac{z}{t}\right)^{\alpha-1} f(t)\, dt, (\alpha > 0, \beta > 0; z \in U^*), \quad (1.4)$$

for $f(z) \in \Sigma$ given by (1.1), we have

$$\mathcal{P}_\beta^\alpha f(z) = z^{-1} + \sum_{k=1}^\infty \left(\frac{\beta}{k+\beta+1}\right)^\alpha a_k z^k, \qquad (\alpha, \beta > 0; z \in U^*). \quad (1.5)$$

From (1.4), we note that

$$z\left(\mathcal{P}_\beta^\alpha f(z)\right)' = \beta \mathcal{P}_\beta^{\alpha-1} f(z) - (\beta+1)\mathcal{P}_\beta^\alpha f(z), \quad (\alpha, \beta > 0; z \in U^*). \quad (1.6)$$

In this paper, we will derive several subordination, superordination and sandwich results involving the operator $\mathcal{P}_\beta^\alpha f(z)$.

1. **Preliminary results.**

**Definition (2.1) [8]:** Let $Q$ be the set of all functions $q$ that are analytic and injective on $\overline{U} \setminus \mathrm{E}(q)$, where

$$\mathrm{E}(q) = \left\{\zeta \in \partial U : \lim_{z \to \zeta} q(z) = \infty\right\} \quad (2.1)$$

and are such that $q'(\zeta) \neq 0$ for $\zeta \in \partial U \setminus E(q)$.

To establish our main results, we need the following lemmas.

**Lemma (2.2) [7]:** Let $q$ be univalent in $U$ and let $\emptyset$ and $\theta$ be analytic in the domain $D$ containing $q(U)$ with $\emptyset(w) \neq 0$, when $w \in q(U)$. Set

$$Q(z) = zq'(z)\emptyset(q(z)) \text{ and } h(z) = \theta(q(z)) + Q(z),$$

suppose that

1 − Q is starlike univalent in U,

$$2 - Re\left(\frac{zh'(z)}{q(z)}\right) = Re\left(\frac{Q'(q(z))}{\emptyset(q(z))} + \frac{zQ'(z)}{Q(z)}\right) > 0, z \in U.$$

If $p$ is analytic in $U$ with $p(0) = q(0)$, $p(U) \subseteq D$ and

$$\theta(p(z)) + zp'(z)\emptyset(p(z)) \prec \theta(q(z)) + zq'(z)\emptyset(q(z)),$$

then $p \prec q$, and $q$ is the best dominant.

**Lemma (2.3) [10]:** Let $q$ be univalent in $U$ and $\psi \in \mathbb{C}$, $\gamma \in \mathbb{C}^*$ with $q(0) = 1$ and

$$Re\left(1 + \frac{zq''(z)}{q'(z)}\right) > \max\left\{0, -Re\left(\frac{\psi}{\gamma}\right)\right\}. \tag{2.2}$$

If $p$ is analytic in U and

$$\psi p(z) + \gamma zp'(z) \prec \psi q(z) + \gamma zq'(z), \tag{2.3}$$

then $p \prec q$ and $q$ is the best dominate.

**Lemma (2.4) [2]:** Let q be convex univalent in the unit disk U and let $\vartheta$ and $\varphi$ be analytic in a domain D containing $q(U)$. Suppose that

$1 - Re\left\{\frac{\vartheta'(q(z))}{\varphi(q(z))}\right\} > 0$ for $z \in U$;

$2 - zq'(z)\varphi(q(z))$ is starlike univalent in $z \in U$.

If $p \in \mathcal{H}[q(0), 1] \cap Q$, with $p(U) \subseteq D$, and $\vartheta(p(z)) + zp'(z)\varphi(q(z))$ is univalent in U, and

$$\vartheta(q(z)) + zq'(z)\varphi(z) \prec \vartheta(p(z)) + zp'(z)\varphi(p(z)), \tag{2.4}$$

then $q \prec p$ and $q$ is the best subordinant.

**Lemma (2.5) [10]:** Let $q$ be convex univalent in U with $q(0) = 1$. Let $\gamma_i \in \mathrm{C}(i = 1,2)$ and $\gamma_2 \neq 0$. Further, assume that $Re\left\{\frac{\gamma_1}{\gamma_2}\right\} > 0$. If $p \in \mathcal{H}[q(0), 1] \cap Q$ and $\gamma_1 p(z) + \gamma p'(z)$ is univalent in U, then

$$\gamma_1 q(z) + \gamma_2 zq'(z) \prec \gamma_1 p(z) + \gamma_2 p'(z), \tag{2.5}$$

which implies that $q \prec p$ and $q$ is the best subordinant.

**Lemma (2.6) [6]:** Let $q$ be univalent in $U$,

$$q(z) = \frac{1+Az}{1-Az}, \text{ with } A \in (-1,0) \cup (0,1), \tag{2.6}$$

and let $\sigma \in (0,1]$, $\alpha, \beta > 0$, such that

$$\frac{2\alpha}{\sigma}\frac{1+A}{1-A} + \frac{\beta}{\sigma}\frac{1+A}{1-A} > 0. \tag{2.7}$$

If $p$ univalent in $U$ with $p(0) = q(0) = 1$, and

$$\alpha p^2(z) + \beta p(z) + \sigma z p'(z) \prec \alpha \left(\frac{1+Az}{1-Az}\right)^2 + \beta \left(\frac{1+Az}{1-Az}\right) + \sigma \frac{2Az}{(1-Az)^2}, \tag{2.8}$$

then $p(z) \prec q(z)$.

### 3. Subordination Results

**Theorem (3.1):** Let $q$ be univalent in the unit disk $U$ with $q(0) = 1$ and $q(z) \neq 0$ for all $z \in U$. Let $\mu, \lambda \in \mathbb{C}^*$, $\gamma \in \mathbb{C}$ and $f \in \Sigma$. Suppose that $\frac{zq'(z)}{q(z)}$ is starlike univalent in $U$ and $f, q$ satisfy the next conditions

$$\frac{(1-\lambda)z\mathcal{P}_\beta^{\alpha-1}f(z) + \lambda z\mathcal{P}_\beta^\alpha f(z)}{\lambda} \neq 0, \tag{3.1}$$

and

$$Re\left(\frac{zq''(z)}{q'(z)} - \frac{zq'(z)}{q(z)} + 1\right) > 0, \tag{3.2}$$

if

$$1 + \gamma\beta\mu \left(\frac{\lambda\mathcal{P}_\beta^\alpha f(z) + (1-2\lambda)\mathcal{P}_\beta^{\alpha-1}f(z) + (\lambda-1)\mathcal{P}_\beta^{\alpha-2}f(z)}{(\lambda-1)\mathcal{P}_\beta^{\alpha-1}f(z) - \lambda\mathcal{P}_\beta^\alpha f(z)}\right) \prec 1 + \gamma\frac{zq'(z)}{q(z)}, \tag{3.3}$$

then

$$\left(\frac{(1-\lambda)z\mathcal{P}_\beta^{\alpha-1}f(z) + \lambda z\mathcal{P}_\beta^\alpha f(z)}{\lambda}\right)^\mu \prec q(z), \tag{3.4}$$

and $q$ is the best dominant.

Proof: We begin by setting

$$\left(\frac{(1-\lambda)z\mathcal{P}_\beta^{\alpha-1}f(z) + \lambda z\mathcal{P}_\beta^\alpha f(z)}{\lambda}\right)^\mu = k(z), \qquad (3.5)$$

then the function $k$ is analytic in $U$ and $q(0) = 1$, and differentiating (3.5) logarithmically with respect to $z$, we get

$$\frac{zk'(z)}{k(z)} = \mu\left(\frac{(1-\lambda)z\left(\mathcal{P}_\beta^{\alpha-1}f\right)'(z) + \lambda z\left(\mathcal{P}_\beta^\alpha f\right)'(z)}{(1-\lambda)\mathcal{P}_\beta^{\alpha-1}f(z) + \lambda\mathcal{P}_\beta^\alpha f(z)} + 1\right). \qquad (3.6)$$

Now, in view of (1.5), we obtain

$$\frac{zk'(z)}{k(z)} = \mu\beta\left(\frac{\lambda\mathcal{P}_\beta^\alpha f(z) + (1-2\lambda)\mathcal{P}_\beta^{\alpha-1}f(z) + (\lambda-1)\mathcal{P}_\beta^{\alpha-2}f(z)}{(\lambda-1)\mathcal{P}_\beta^{\alpha-1}f(z) - \lambda\mathcal{P}_\beta^\alpha f(z)}\right), \qquad (3.7)$$

by setting

$$\theta(w) = 1 \text{ and } \emptyset(w) = \frac{\gamma}{w},$$

it can easily observed that $\theta(w)$ in analytic in $\mathbb{C}$ and $\emptyset(w) \neq 0$ is an analytic in $\mathbb{C}^*$.

Moreover, we let

$$Q(z) = zq'(z)\emptyset(q(z)) = \gamma\frac{zq'(z)}{q(z)}, \qquad (3.8)$$

and

$$h(z) = \theta(q(z)) + Q(z) = 1 + \gamma\frac{zq'(z)}{q(z)}. \qquad (3.9)$$

We find that $Q(z)$ is starlike univalent in $U$, and from (3.2)

$$Re\left(\frac{zh'(z)}{Q(z)}\right) = Re\left(\frac{zq''(z)}{q'(z)} - \frac{zq'(z)}{q(z)} + 1\right) > 0, \qquad (3.10)$$

and by using Lemma (2.2), we deduce that the subordination (3.3) implies $k(z) \prec q(z)$, and the function $q$ is the best dominant of (3.3).

Upon setting $q(z) = e^{\tau z}, |\tau| \leq 1$ in Theorem (3.1) we get the following result.

**Corollary (3.2):** Let $f \in \Sigma$ and $|\tau| \leq 1$, also the condition (3.2) is satisfied. If

$$1 + \gamma\beta\mu\left(\frac{\lambda\mathcal{P}_\beta^\alpha f(z) + (1-2\lambda)\mathcal{P}_\beta^{\alpha-1} f(z) + (\lambda-1)\mathcal{P}_\beta^{\alpha-2} f(z)}{(\lambda-1)\mathcal{P}_\beta^{\alpha-1} f(z) - \lambda\mathcal{P}_\beta^\alpha f(z)}\right) \prec 1 + \tau\gamma z, \qquad (3.11)$$

then

$$\left(\frac{(1-\lambda)z\mathcal{P}_\beta^{\alpha-1} f(z) + \lambda z\mathcal{P}_\beta^\alpha f(z)}{\lambda}\right)^\mu \prec e^{\tau z}, \qquad (3.12)$$

and $e^{\tau z}$ is the best dominant.

Hence, for the particular case $\tau = \lambda = 1$, we have the following result.

**Corollary (3.3):** let $f \in \Sigma$ satisfies the subordination

$$1 - \gamma\mu\beta\left(1 - \frac{\mathcal{P}_\beta^{\alpha-1} f(z)}{\mathcal{P}_\beta^\alpha f(z)}\right) \prec 1 + \gamma z, \qquad (3.13)$$

then

$$\left(z\mathcal{P}_\beta^\alpha f(z)\right)^\mu \prec e^z. \qquad (3.14)$$

Putting $q(z) = \left(\frac{1+z}{1-z}\right)^\rho, 0 < \rho \leq 1 \text{ and } \lambda = 1$ in Theorem (3.1), we obtain the following result.

**Corollary (3.4):** let $0 < \rho \leq 1$ and $f \in \Sigma$ satisfies the subordination

$$1 - \gamma\mu\beta\left(1 - \frac{\mathcal{P}_\beta^{\alpha-1} f(z)}{\mathcal{P}_\beta^\alpha f(z)}\right) \prec 1 + 2\rho\gamma \frac{z}{1-z^2}, \qquad (3.15)$$

then

$$\left(z\mathcal{P}_\beta^\alpha f(z)\right)^\mu \prec \left(\frac{1+z}{1-z}\right)^\rho \qquad (3.16)$$

and $\left(\frac{1+z}{1-z}\right)^\rho$ is the best dominant.

**Theorem (3.5):** Let $q$ be convex univalent in the unit disk U with $q(0) = 1$, let $\gamma > 0, \mu \in \mathbb{C}^*, \eta, \delta \in \mathbb{C}, f \in \Sigma$ and suppose that $f$ and $q$ satisfy the following conditions

$$z\mathcal{P}_\beta^\alpha f(z) \neq 0, \tag{3.17}$$

and

$$Re\left(\frac{\eta}{\gamma} + \frac{2\delta}{\gamma} q(z)\right) > 0. \tag{3.18}$$

If

$$\psi(z) \prec \delta q^2(z) + \eta q(z) + \gamma z q'(z), \tag{3.19}$$

where

$$\psi(z) = \left(z\mathcal{P}_\beta^\alpha f(z)\right)^\mu \left(\delta \left(z\mathcal{P}_\beta^\alpha f(z)\right)^\mu + \gamma\mu\beta \left(\frac{\mathcal{P}_\beta^{\alpha-1} f(z)}{\mathcal{P}_\beta^\alpha f(z)} - 1\right) + \eta\right), \tag{3.20}$$

then

$$\left(z\mathcal{P}_\beta^\alpha f(z)\right)^\mu \prec q(z), \tag{3.21}$$

and $q$ is the best dominant.

**Proof:** Let

$$p(z) = \left(z\mathcal{P}_\beta^\alpha f(z)\right)^\mu, \quad z \in U. \tag{3.22}$$

According to (3.17) the function $p(z)$ is analytic in $U$ with $p(0) = 1$. A simple computation shows that

$$\psi(z) = \left(z\mathcal{P}_\beta^\alpha f(z)\right)^\mu \left(\delta \left(z\mathcal{P}_\beta^\alpha f(z)\right)^\mu + \gamma\mu\beta \left(\frac{\mathcal{P}_\beta^{\alpha-1} f(z)}{\mathcal{P}_\beta^\alpha f(z)} - 1\right) + \eta\right)$$

$$= \delta p^2(z) + \eta p(z) + \gamma z p'(z), \tag{3.23}$$

to prove our result use Lemma (2.2), consider in this lemma $\theta(w) = \delta w^2 + \eta w$ and $\phi(w) = \gamma$, then $\theta$ is analytic in $\mathbb{C}$ and $\phi$ is analytic in $\mathbb{C}^*$. Also, if we let

$$Q(z) = zq'(z)\phi(z) = \gamma z q'(z), \tag{3.24}$$

and

$$h(z) = \theta(q(z)) + Q(z) = \delta q^2(z) + \eta q(z) + \gamma z q'(z). \tag{3.25}$$

Then, the assumption $q$ is convex would yield $Q$ is a starlike function in $U$. From (3.18) we have

$$Re\left(\frac{zh'(z)}{Q(z)}\right) = Re\left(1 + \frac{\eta}{\gamma} + \frac{2\delta}{\gamma}q(z) + \gamma\frac{q''(z)}{q'(z)}\right) > 0, \tag{3.26}$$

and by using Lemma (2.2), we deduce that the subordination (3.19) implies that

$p(z) \prec q(z)$, and the function $q$ is the best dominant.

Taking $q(z) = \frac{1+Az}{1-Az}, A \in (-1,0) \cup (0,1)$ in Theorem (35) and using Lemma (2.6), we have the next result.

**Corollary (3.6):** Let $q(z) = \frac{1+Az}{1-Az}$ with $A \in (-1,0) \cup (0,1), \eta, \delta > 0, \gamma \in (0,1]$, such that

$$Re\left(\frac{\eta}{\gamma} + \frac{2\delta}{\gamma}\frac{1+Az}{1-Az} + \frac{1+Az}{1-Az}\right) > 0. \tag{3.27}$$

If $f \in \Sigma$ satisfies the subordination

$$\psi(z) \prec \delta\left(\frac{1+Az}{1-Az}\right)^2 + \eta\frac{1+Az}{1-Az} + \gamma\frac{2Az}{(1-Az)^2}, \tag{3.28}$$

where

$$\psi(z) = \left(z\mathcal{P}_\beta^\alpha f(z)\right)^\mu \left(\delta\left(z\mathcal{P}_\beta^\alpha f(z)\right)^\mu + \gamma\mu\beta\left(\frac{\mathcal{P}_\beta^{\alpha-1}f(z)}{\mathcal{P}_\beta^\alpha f(z)} - 1\right) + \eta\right), \tag{3.29}$$

then

$$\left(z\mathcal{P}_\beta^\alpha f(z)\right)^\mu \prec \left(\frac{1+Az}{1-Az}\right). \tag{3.30}$$

Putting $q(z) = e^{\tau z}, |\tau| \leq 1$ in Theorem (3.5), we get the following result.

**Corollary (3.7):** Assume that

$$Re\left(1+\frac{\eta}{\gamma}+z\tau+\frac{2\delta}{\gamma}e^{\tau z}\right)>0, \qquad \gamma>0 \qquad (3.31)$$

if $f \in \Sigma$ is satisfies the subordination

$$\left(z\mathcal{P}_\beta^\alpha f(z)\right)^\mu\left(\delta\left(z\mathcal{P}_\beta^\alpha f(z)\right)^\mu+\gamma\mu\beta\left(\frac{\mathcal{P}_\beta^{\alpha-1}f(z)}{\mathcal{P}_\beta^\alpha f(z)}-1\right)+\eta\right) \prec e^{\tau z}(\delta e^{\tau z}+\eta+\gamma\tau z), \quad (3.32)$$

then

$$\left(z\mathcal{P}_\beta^\alpha f(z)\right)^\mu \prec e^{\tau z}, \qquad (3.33)$$

and $e^{\tau z}$ is the best dominant.

**Theorem (3.8):** Let $q(z)$ be univalent in $U$ with $q(0)=1$. Suppose that

$$Re\left(1+\frac{zq''(z)}{q'(z)}\right) > \max\left\{0, Re\left(\frac{m\beta}{\ell}\right)\right\}, m \in \mathbb{C}, (\beta, \ell \in \mathbb{C}^*), z \in U, \qquad (3.34)$$

if $f \in \Sigma$ is satisfies the subordination

$$\Phi(z) \prec mq(z)-\frac{\ell}{\beta}zq'(z), \qquad (3.35)$$

where

$$\Phi(z)=\left(\frac{(1-\lambda)z\mathcal{P}_\beta^{\alpha-1}f(z)+\lambda z\mathcal{P}_\beta^\alpha f(z)}{\lambda}\right)^\mu\Bigg[m$$
$$+\mu\ell\left(\frac{\lambda\mathcal{P}_\beta^\alpha f(z)+(1-2\lambda)\mathcal{P}_\beta^{\alpha-1}f(z)+(\lambda-1)\mathcal{P}_\beta^{\alpha-2}f(z)}{(1-\lambda)\mathcal{P}_\beta^{\alpha-1}f(z)+\lambda\mathcal{P}_\beta^\alpha f(z)}\right)\Bigg], \qquad (3.36)$$

then

$$\left(\frac{(1-\lambda)z\mathcal{P}_\beta^{\alpha-1}f(z)+\lambda z\mathcal{P}_\beta^\alpha f(z)}{\lambda}\right)^\mu \prec q(z), \qquad (3.37)$$

and $q$ is the best dominant.

**Proof:** Let $p(z)$ be defined by (3.5). Then simple computations show that

$$\left(\frac{(1-\lambda)z\mathcal{P}_\beta^{\alpha-1}f(z) + \lambda z\mathcal{P}_\beta^\alpha f(z)}{\lambda}\right)^\mu \left[m + \mu\ell\left(\frac{\lambda\mathcal{P}_\beta^\alpha f(z) + (1-2\lambda)\mathcal{P}_\beta^{\alpha-1}f(z) + (\lambda-1)\mathcal{P}_\beta^{\alpha-2}f(z)}{(1-\lambda)\mathcal{P}_\beta^{\alpha-1}f(z) + \lambda\mathcal{P}_\beta^\alpha f(z)}\right)\right]$$

$$= mp(z) - \frac{\ell}{\beta}zp'(z). \tag{3.38}$$

Thus the subordination (3.35) is equivalent to

$$mp(z) - \frac{\ell}{\beta}zp'(z) \prec mq(z) - \frac{\ell}{\beta}zq'(z). \tag{3.39}$$

Applying Lemma (2.3) with $\psi = m, \gamma = \frac{-\ell}{\beta}$. The proof of theorem is complete.

Taking $q(z) = \left(\frac{1+z}{1-z}\right)^\rho, 0 < \rho \leq 1$ and $\lambda = 1$ in Theorem (3.8), we obtain the following result.

**Corollary (3.9):** Let $0 < \rho \leq 1$ and suppose that

$$Re\left(\frac{1 + 2\rho z + z^2}{1 - z^2}\right) > \max\left\{0, Re\left(\frac{m\beta}{\ell}\right)\right\}, \tag{3.40}$$

if $f \in \Sigma$ is satisfies the subordination

$$\left(z\mathcal{P}_\beta^\alpha f(z)\right)^\mu \left(m + \ell\mu\left(1 + \frac{\mathcal{P}_\beta^{\alpha-1}f(z)}{\mathcal{P}_\beta^\alpha f(z)}\right)\right) \prec \left(\frac{1+z}{1-z}\right)^\rho \left(m - \frac{\ell}{\mathcal{B}}\left(\frac{2\rho z}{1-z^2}\right)\right), \tag{3.41}$$

then

$$\left(z\mathcal{P}_\beta^\alpha f(z)\right)^\mu \prec \left(\frac{1+z}{1-z}\right)^\rho, \tag{3.42}$$

and $\left(\frac{1+z}{1-z}\right)^\rho$ is the best dominant.

Putting, $q(z) = \frac{1+Az}{1+Bz}, (-1 \leq B < A \leq 1)$ and $\lambda = 1$ in Theorem (3.8), the condition (3.34) reduces to

$$Re\left(\frac{1 - Bz}{1 + BA}\right) > \max\left\{0, Re\left(\frac{m\beta}{\ell}\right)\right\}, \beta \neq 0. \tag{3.43}$$

It is easy to verify that the function $\varphi(\delta) = \frac{1-\delta}{1+\delta}, |\delta| < |\beta|$, is convex in $U$ and since $\varphi(\bar{\delta}) = \overline{\varphi(\delta)}$ for all $|\delta| < |\beta|$, it follows that $\varphi(\delta)$ is a convex domain symmetric with respect to real axis, hence

$$\inf\left\{Re\left(\frac{1-Bz}{1+Bz}\right), z \in U\right\} = \frac{1-|B|}{1+|B|} > 0. \tag{3.44}$$

We obtain the following corollary.

**Corollary (3.10):** Let $(-1 \leq B < A \leq 1)$ and $\frac{1-|B|}{1+|B|} \geq \max\left\{0, Re\left(\frac{m\beta}{\ell}\right)\right\}$,

if $f \in \Sigma$ is satisfies the subordination

$$\left(z\mathcal{P}_\beta^\alpha f(z)\right)^\mu \left(m + \ell\mu\left(1 - \frac{\mathcal{P}_\beta^{\alpha-1} f(z)}{\mathcal{P}_\beta^\alpha f(z)}\right)\right) \prec m\frac{1+Az}{1+Bz} - \frac{(A-B)z}{(1+Bz)^2}, \tag{3.45}$$

then

$$\left(z\mathcal{P}_\beta^\alpha f(z)\right)^\mu \prec \frac{1+Az}{1+Bz}, \tag{3.46}$$

and $\frac{1+Az}{1+Bz}$ is the best dominant.

## 4. Superordination Results

**Theorem (4.1)** Let q be convex univalent in U with $q(0) = 1. Let\ \gamma > 0, \mu \in \mathbb{C}^*, \eta, \delta \in \mathbb{C}$ and $f \in \Sigma$ .Suppose that

$$Re\left\{\left(\frac{\eta}{\gamma} + \frac{2\delta}{\gamma}q(z)\right)q'(z)\right\} > 0 \tag{4.1}$$

and $f$ satisfies the next conditions

$$z\mathcal{P}_\beta^\alpha f(z) \neq 0, \qquad (\alpha, \beta > 0; z \in U) \tag{4.2}$$

and

$$\left(z\mathcal{P}_\beta^\alpha f(z)\right)^\mu \in H[q(0), 1] \cap Q, \tag{4.3}$$

also, if the function $\psi(z)$ defined by (3.20) is univalent in U and the following superordination condition

$$\delta q^2(z) + \eta q(z) + \gamma z q'(z) \prec \psi(z), \tag{4.4}$$

holds, then

$$q(z) \prec \left(z\mathcal{P}_\beta^\alpha f(z)\right)^\mu, \tag{4.5}$$

and $q$ is the best subordinant.

**Proof**: Let

$$g(z) = \left(z\mathcal{P}_\beta^\alpha f(z)\right)^\mu, \tag{4.6}$$

then, after computation, we get

$$\delta g^2(z) + \eta g(z) + \gamma z g'(z) = \psi, \tag{4.7}$$

this implies

$$\delta q^2(z) + \eta q(z) + \gamma z q'(z) \prec \delta g^2(z) + \eta g(z) + \gamma z g'(z). \tag{4.8}$$

by setting

$$\vartheta(w) = \delta w^2 + \eta w \quad \text{and} \quad \varphi(w) = \gamma,$$

it can easily observed that $\vartheta(w)$ is analytic in $\mathbb{C}$, and $\varphi(w) \neq 0$ is an analytic in $\mathbb{C}^*$.

Also, we obtain

$$Re\left\{\frac{\vartheta'(q(z))}{\varphi(z)}\right\} = Re\left\{\left(\frac{\eta}{\gamma} + \frac{2\delta}{\gamma}q(z)\right)q'(z)\right\} > 0. \tag{4.9}$$

Therefore, by Lemma (2.4), we have

$$q(z) \prec \left(z\mathcal{P}_\beta^\alpha f(z)\right)^\mu. \tag{4.10}$$

Taking $q(z) = e^{\tau z}, |\tau| \leq 1$ in Theorem (4.1), we have the following corollary.

**Corollary (4.2)**: Let $Re\left\{\left(\frac{\eta}{\gamma} + \frac{2\delta}{\gamma}e^{\tau z}\right)\tau e^{\tau z}\right\} > 0$ and $f \in \Sigma$ such that $\left(z\mathcal{P}_\beta^\alpha f(z)\right)^\mu \in H[q(0), 1] \cap Q$. If the function $\psi(z)$ defined by (3.20) is univalent in U and satisfied the following superordination condition

$$e^{\tau z}(\delta e^{\tau z} + \eta + \gamma \tau z) \prec \left(z\mathcal{P}_\beta^\alpha f(z)\right)^\mu \left(\delta \left(z\mathcal{P}_\beta^\alpha f(z)\right)^\mu + \gamma \mu \beta \left(\frac{\mathcal{P}_\beta^{\alpha-1} f(z)}{\mathcal{P}_\beta^\alpha f(z)} - 1\right) + \eta\right), \tag{4.11}$$

then

$$e^{\tau z} \prec \left(z\mathcal{P}_{\beta}^{\alpha}f(z)\right)^{\mu}, \tag{4.12}$$

and $e^{\tau z}$ is best subordinant.

Now, by appealing to Lemma (2.5) it easily to prove the following theorem.

**Theorem (4.3)** Let q be convex univalent in U with $q(0) = 1$ and $Re\left(\frac{m\beta}{\ell}\right) < 0, m \in \mathbb{C}, \beta, \ell \in \mathbb{C}^*$, if $f \in \Sigma$ such that

$$\frac{(1-\lambda)z\mathcal{P}_{\beta}^{\alpha-1}f(z) + \lambda z\mathcal{P}_{\beta}^{\alpha}f(z)}{\lambda} \neq 0, \tag{4.13}$$

and

$$\left(\frac{(1-\lambda)z\mathcal{P}_{\beta}^{\alpha-1}f(z) + \lambda z\mathcal{P}_{\beta}^{\alpha}f(z)}{\lambda}\right)^{\mu} \in H[q(0),1] \cap Q. \tag{4.14}$$

If the function $\Phi(z)$ defined by (3.36) is univalent in U and the following superordination condition

$$mq(z) - \frac{\ell}{\beta}zq'(z) \prec \Phi(z), \tag{4.15}$$

holds, then

$$q(z) \prec \left(\frac{(1-\lambda)z\mathcal{P}_{\beta}^{\alpha-1}f(z) + \lambda z\mathcal{P}_{\beta}^{\alpha}f(z)}{\lambda}\right)^{\mu}, \tag{4.16}$$

and $q(z)$ is the best subordinant.

Putting, $q(z) = \frac{1+Az}{1+Bz}, (-1 \leq B < A \leq 1)$ and $\lambda = 1$ in Theorem (4.3), we get the following result.

**Corollary (4.4)**: Let $Re\left(\frac{m\beta}{\ell}\right) < 0$ and $f \in \Sigma$ such that $\left(z\mathcal{P}_{\beta}^{\alpha}f(z)\right)^{\mu} \in \mathcal{H}[q(0),1] \cap Q$, $\left(z\mathcal{P}_{\beta}^{\alpha}f(z)\right)^{\mu}\left(m + \ell\mu\left(1 + \frac{\mathcal{P}_{\beta}^{\alpha-1}f(z)}{\mathcal{P}_{\beta}^{\alpha}f(z)}\right)\right)$ is univalent in U and satisfied the following superordination condition

$$m\frac{1+Az}{1+Bz} - \frac{(A-B)z}{(1+Bz)^2} \prec \left(z\mathcal{P}_\beta^\alpha f(z)\right)^\mu \left(m + \ell\mu\left(1 + \frac{\mathcal{P}_\beta^{\alpha-1}f(z)}{\mathcal{P}_\beta^\alpha f(z)}\right)\right), \qquad (4.17)$$

then

$$\frac{1+Az}{1+Bz} \prec \left(z\mathcal{P}_\beta^\alpha f(z)\right)^\mu,$$

and $\frac{1+Az}{1+Bz}$ is the best subordinant.

Taking $q(z) = \left(\frac{1+z}{1-z}\right)^\rho, 0 < \rho \leq 1$ and $\lambda = 1$, in Theorem (4.3), we obtain the following result.

**Corollary (4.5) :** Let $Re\left(\frac{m\beta}{\ell}\right) < 0$ and $f \in \Sigma$ such that $\left(z\mathcal{P}_\beta^\alpha f(z)\right)^\mu \in H[q(0), 1] \cap Q$,

$\left(z\mathcal{P}_\beta^\alpha f(z)\right)^\mu \left(m + \ell\mu\left(1 - \frac{\mathcal{P}_\beta^{\alpha-1}f(z)}{\mathcal{P}_\beta^\alpha f(z)}\right)\right)$ is univalent in U and satisfied the following superordination condition

$$\left(\frac{1+z}{1-z}\right)^\rho \left(m - \frac{\ell}{\mathcal{B}}\left(\frac{2\rho z}{1-z^2}\right)\right) \prec \left(z\mathcal{P}_\beta^\alpha f(z)\right)^\mu \left(m + \ell\mu\left(1 - \frac{\mathcal{P}_\beta^{\alpha-1}f(z)}{\mathcal{P}_\beta^\alpha f(z)}\right)\right), \qquad (4.18)$$

then

$$\left(\frac{1+z}{1-z}\right)^\rho \prec \left(z\mathcal{P}_\beta^\alpha f(z)\right)^\mu, \qquad (4.19)$$

and $\left(\frac{1+z}{1-z}\right)^\rho$ is the best subordinant.

## 5. Sandwich Result

Combining Theorem (3.5) with Theorem (4.1) and Theorem (3.8) with Theorem (4.3), we obtain, respectively, the following two sandwich results:

**Theorem (5.1)** Let $q_1$ be convex univalent in U with $q_1(0) = 1$, and satisfies (4.1) and let $q_2$ be univalent in U with $q_2(0) = 1$ and satisfies (3.18). Further suppose that $\gamma > 0, \mu \in \mathbb{C}^*, \eta, \delta \in \mathbb{C}$,

$$z\mathcal{P}_\beta^\alpha f(z) \neq 0, \qquad (5.1)$$

and

$$\left(z\mathcal{P}_\beta^\alpha f(z)\right)^\mu \in H[q(0), 1] \cap Q. \tag{5.2}$$

If the function $\psi(z)$ defined by (3.20) is univalent in U and

$$\delta q_1^2(z) + \eta q_1(z) + \gamma z q'_1(z) \prec \psi(z) \prec \delta q_2^2(z) + \eta q_2(z) + \gamma z q'_2(z), \tag{5.3}$$

then

$$q_1(z) \prec \left(z\mathcal{P}_\beta^\alpha f(z)\right)^\mu \prec q_2(z), \tag{5.4}$$

and $q_1$ and $q_2$ are respectively, the best subordinant and best dominant.

**Theorem (5.2)** Let $q_1$ be convex univalent in U with $q_1(0) = 1$, and let $q_2$ be univalent in U. Suppose that, $Re\left(\frac{m\beta}{\ell}\right) < 0, m \in \mathbb{C}, \beta, \ell \in \mathbb{C}^*$, and $q_2$ satisfies (3.34), let

$$\left(\frac{(1-\lambda)z\mathcal{P}_\beta^{\alpha-1} f(z) + \lambda z\mathcal{P}_\beta^\alpha f(z)}{\lambda}\right)^\mu \in H[q(0), 1] \cap Q, \tag{5..5}$$

and the function $\Phi(z)$ defined by (3.36) is univalent in U, if if $f \in \Sigma$ satisfies

$$mq_1(z) - \frac{\ell}{\beta} z q_1'(z) \prec \Phi(z) \prec mq_2(z) - \frac{\ell}{\beta} z q_2'(z), \tag{5.6}$$

then

$$q_1(z) \prec \left(\frac{(1-\lambda)z\mathcal{P}_\beta^{\alpha-1} f(z) + \lambda z\mathcal{P}_\beta^\alpha f(z)}{\lambda}\right)^\mu \prec q_2(z), \tag{5.7}$$

and $q_1$ and $q_2$ are the best subordinant and dominant, respectively.